%% file: draft.tex
\documentclass{amsart}

\usepackage{header}


\title{Deep Learning of Topological Equivalence Classes}

\author[Linden Disney-Hogg]{Linden Disney-Hogg$^\ast$ \orcidB{}}
\address{
Department of Mathematics \\ Imperial College London \\ 
London SW7 2AZ, U.K.
}
\email{a.disneyhogg@imperial.ac.uk}

\author[Max Zuev]{Max Zuev \orcidA{}}
\address{
Mathematical Institute \\ The University of Oxford \\ 
Oxford OX2 6GG, U.K.
}
\email{maxim.zuev@maths.ox.ac.uk}

\thanks{$^\ast$ Corresponding author}

\thanks{\textbf{Acknowledgements.} MZ was supported by a University of Leeds School of Mathematics Undergraduate Summer Research Bursary. The research of LDH was supported by a UK Engineering and Physical Sciences Research Council (EPSRC) doctoral prize fellowship, and a Chapman-Schmidt AI in Science postdoctoral fellowship, a Schmidt Sciences program. Computing for this work was undertaken on the Aire HPC system at the University of Leeds, UK.}

\thanks{\textbf{Data Availability.} The datasets generated during the current study and the code for their creation/analysis are openly available at \url{https://github.com/Max-Zuev/PermutationGroupDeepLearning}.}

\thanks{\textbf{Statements and Declarations.} The authors have no relevant financial or non-financial interests to disclose. \\
The generative AI tool Claude Code was used in the implementation and optimisation of the code written for this project. All code was subsequently checked by humans.}

\date{October 2026}

\begin{document}

\begin{abstract}
We demonstrate an application of deep learning to the problem of identifying whether two generating vectors correspond to topologically equivalent finite group actions on Riemann surfaces, showing how architectures developed for semantic understanding have use for pure mathematical problems. In doing so we discuss how desired algebraic structure influences model design, and explore which structural properties of topological equivalence classes the model has successfully learnt. 
\end{abstract}

\maketitle

\section{Introduction}

A fundamental problem in the study of Riemann surfaces is an understanding of their moduli, and one key approach to this is via describing the possible holomorphic actions of a finite group $G$ on a smooth compact connected Riemann surface $S$, defined as monomorphisms $\epsilon : G \to \Aut(S)$. A foundational result in this direction is the theorem of Hurwitz showing that when $S$ has genus $g \geq 2$, $\abs{\Aut(S)} \leq 84(g-1)$ \cite{Hurwitz1892}.

Associated to the image of a group action $\mathcal{G} = \epsilon(G) \leq \Aut(S)$ one has the quotient map $\pi : S \to S/\mathcal{G}$. 
\begin{definition}[\cite{Lange2022}]
    Given an action of $G$ on $S$, define the \textbf{quotient genus} to be $g_0 = g(S/\mathcal{G}) \geq 0$. Moreover, letting $Q_i$, $i=1, \dots, r$, denote the branch points of $\pi$, for any $P_i \in \pi^{-1}(Q_i)$ denote with $c_i = \abs{\mathcal{G}_{P_i}} >1$ the order of the stabiliser of $P_i$ in $\mathcal{G}$. We call the integer tuple $(g_0; c_1, \dots, c_r)$ the \textbf{signature} of the action. 

    Moreover, given any finite group $G$ and signature $(g_0; c_1, \dots, c_r)$, define a \textbf{generating vector} associated to the signature to be $\pbrace{\alpha_1, \beta_1, \dots, \alpha_{g_0}, \beta_{g_0}, \gamma_1, \dots, \gamma_r} \subset G$ such that 
    \begin{itemize}
        \item $\lbrace\alpha_1, \beta_1, \dots, \alpha_{g_0}, \beta_{g_0}, \gamma_1, \dots, \gamma_r\rbrace$ generates $G$,
        \item $\gamma_i$ has order $c_i$, and 
        \item $\prod_{i=1}^{g_0} (\alpha_i \beta_i \alpha_i^{-1} \beta_i^{-1}) \prod_{j=1}^r \gamma_j = 1$. 
    \end{itemize}
\end{definition}
The Riemann existence theorem states that, given a finite group $G$, a signature $(g_0; c_1, \dots, c_r)$, and a generating vector in $G$ associated to this signature, then there exists an action of $G$ on a surface $S$ of genus $g$ if and only if the Riemann-Hurwitz equation
\[
2g - 2 = \abs{G} \psquare{(2g_0 - 2) + \sum_{i=1}^r \pround{1- \frac{1}{c_i}}}
\]
holds \cite{Broughton1990}. In this work we shall restrict ourselves to investigating group actions where $g_0=0$ (and hence $r \geq 3$), but this is not a theoretical barrier to the application of our methods. 

\begin{remark}
	Given the data of the character table of $G$ one may use character theory to estimate the difficulty of finding all generating vectors of $G$ for a particular signature, and the expected number. Denote with $C_i$ a conjugacy class in $G$, let $x_i \in C_i$, and denote with $\lbrace \chi_i \rbrace$ the irreducible characters of $G$. It is then known the size of the conjugacy classes may be determined by the orthogonality relation 
	$$
	\frac{1}{|G|}\sum_i \chi_i(x_r) \chi_i(x_s) = \frac{\delta_{rs}}{|C_r|}.
	$$
	Fixing conjugacy classes $C_1, \dots, C_r$ one may imagine an algorithm which searches for generating vectors with $\gamma_i \in C_i$ by looping over the cartesian product $C_1 \times \dots \times C_{r-1}$ and checking whether the inverse of the product is in $C_r$. This is a search over a space of size $|C_1|\times \dots \times |C_{r-1}|$ which may be determined, and the number of solutions is given by Frobenius' formula \cite[Theorem 14.1]{Breuer2000}
	\[
	|\lbrace (\gamma_1, \dots, \gamma_r) \, | \, \gamma_i \in C_i, \, \prod_{i=1}^r \gamma_i = 1 \rbrace | = \frac{\prod_{i=1}^r |C_i|}{|G|} \sum_i \frac{\prod_{j=1}^r \chi_i(x_j)}{\chi_i(1)^{r-2}} .
	\]
\end{remark}

As is natural in mathematics one wants to describe all such group actions up to a concept of equivalence; many such concepts of equivalence have been considered \cite{Broughton2022}, but in this paper we shall focus on perhaps the most well studied, namely \textit{topological equivalence}.

\begin{definition}[\cite{Kuribayashi1987}, \cite{Broughton2022}]
    Actions $\epsilon_i$ of $G$ on Riemann surfaces $S_i$ ($i=1,2$) of the same genus are said to be \textbf{topologically equivalent} if there exists $h : S_1 \to S_2$ an orientation preserving homeomorphism and $\sigma \in \Aut(G)$ such that $\forall g \in G$, $\epsilon_2(g) = h \circ \epsilon_1(\sigma(g)) \circ h^{-1}$. 
\end{definition}

This definition of equivalence is particularly important because topological equivalence classes biject with conjugacy classes of finite subgroups of the mapping class group, a fundamental object in many areas of mathematics \cite{Broughton1990}. It is clear that the signature of a group action is an invariant of its topological equivalence class up to permutation of the $c_i$. Through the theory of surface kernel epimorphisms one can show that topological equivalence classes correspond to orbits of generating vectors under the action \cite{Broughton1990, Behn2023} of
\begin{itemize}
    \item $\Aut(G)$, with $\sigma \in \Aut(G)$ acting on $(\gamma_1, \dots, \gamma_r)$ as 
    \[
    \sigma \cdot (\gamma_1, \dots, \gamma_r) = (\sigma(\gamma_1), \dots, \sigma(\gamma_r)),
    \]
    and 
    \item the subgroup of the braid group $\mathbb{B}_r$ generated by $\Phi_i$, $i=1, \dots, r-1$, acting as 
    \[
    \Phi_i \cdot (\gamma_1, \dots, \gamma_r) = (\gamma_1, \dots, \gamma_{i-1}, \gamma_i \gamma_{i+1} \gamma_i^{-1}, \gamma_i, \dots, \gamma_r).
    \]
\end{itemize}
We shall further fix the ordering of the signature $(c_1, \dots, c_r)$, and consider only the subgroup of $\mathbb{B}_r$ which does not permute the signature. 

\begin{example}\label{ex: cyclic gvs}
Consider the example of $G = \langle \gamma \, | \, \gamma^4=1 \rangle \cong C_4$ acting with signature $(0; 4, 4, 4, 4)$ in genus $g=3$. Any generating vector must be of the form 
\[
v = (\gamma^{a_1}, \gamma^{a_2}, \gamma^{a_3}, \gamma^{a_4}), 
\]
where $a_i = 1, 3$ and $\sum_i a_i \equiv 0 \mod 4$. $\Aut(G) = \langle \gamma \mapsto \gamma^3 \rangle$ and any permutation preserves the signature, so we have two equivalence classes of generating vector
\begin{itemize}
    \item two of the $a_i$ equal 1, and two of the $a_i$ equal 3,
    \item all of the $a_i$ are equal. 
\end{itemize}
\end{example}

\begin{example}\label{ex: bad gv}
    Consider the group $G = A_4$ the alternating group on 4 elements acting with signature $(0; 3, 3, 3, 3)$ in genus $g=5$. Any generating vector must consist of four $3$-cycles, two examples are
    \begin{align*}
    v_1 &= ((1,2,4), (1,3,4), (2,3,4), (1,2,3)), \\
    v_2 &=((2,3,4), (1,3,4), (1,2,4), (1,2,3)).
    \end{align*}
    One may check that these are in distinct equivalence classes despite differing by a permutation. 
\end{example}

Given a generating vector it is easy to compute its associated signature but this is not a complete invariant of the equivalence class, that is there exist topologically inequivalent generating vectors with the same signature (as demonstrated in Examples \ref{ex: cyclic gvs} and \ref{ex: bad gv}). There are more refined invariants, such as the character of the induced representation of $G$ on the space $\mathcal{H}^q(S)$ of holomorphic $q$-differentials on $S$, which are similarly quick to compute but are still not complete. The exception to this is in the case where $G$ is cyclic and the `total rotation data system' provides a complete invariant \cite{Broughton2022}.

In contrast, one may compute the entire orbit of a generating vector using modern computer algebra \cite{Behn2023} and use this to determine if two generating vectors are topologically equivalent, but this process can become slow and memory intensive dependent on the size of $\Aut(G) \times \mathbb{B}_r$. Indeed, in general the orbit problem is at least as hard as the graph isomorphism problem \cite{Clarke1998}.

\begin{example}
    For a finite abelian group $G$ it is known that $|\mathrm{Aut}(G)| \geq \phi(|G|)$, where $\phi$ is Euler's totient function. Moreover when $G$ is abelian the action of the braid generator $\Phi_i$ on a generating vector $(\gamma_1, \dots, \gamma_r)$ is simply to swap $\gamma_i$ and $\gamma_{i+1}$, and so the braid action is just that of the symmetric group $S_r$. As such when considering the topological equivalence classes for generating vectors of length $r$ for an abelian group $G$, orbits are generated by the action of a group of at least order $\phi(|G|)r!$.  
\end{example}

Faced with these difficulties we sought to investigate what alternative computational tools could provide insight into the problem of determining whether two generating vectors correspond to topologically equivalent actions. Motivated by recent successes in the application of machine learning to understanding group actions on Riemann surfaces \cite{Braden2025} as well as in using deep learning to understand structure of the symmetric group \cite{Petschack2025} we decided to use a deep learning approach. In particular, we trained a twin model originally designed to decide if two sequences of text share the same `sentiment' \cite{Yang2019}: in our context the generating vectors gave rise to sequences of tokens for which sharing sentiment corresponded to being in the same equivalence class. In doing so we add to a growing corpus of literature exploring the applicability of machine/deep learning techniques to pure mathematics, see \cite{Douglas2025} for a review. The field of Riemann surfaces and particularly those with automorphisms is well placed to benefit from the development of this technology thanks to the wealth of precomputed data available from the LMFDB \cite{lmfdb, Paulhus2023}. 

The final best-performing model successfully compared generating vectors in the test dataset with an accuracy of $85.8\%$. 
This paper provides further evidence of how modern computational developments can be leveraged to gain insight into mathematical problems. We also hope it will serve as an impetus to further explore the interplay between such fields. \S2 provides specific details of the model architecture employed, the training data used, how results were analysed, and other key technical information, while \S3 provides the results of investigations. Finally, \S4 discusses the interpretation of these results, as well as directions of future inquiry.

\section{Methods}

\begin{figure}
    \centering  \input{twin_network_tikz_paper_v2}
    \caption{The twin network architecture. The model has two identical halves, each half processing one generating vector. A tokenised element of a generating vector is passed through three identical processing blocks (with distinct weights) in series: the output of a convolutional encoder is concatenated with the block input, an alignment layer relates the encodings to those of the corresponding generating vector in the opposite half, and the outputs are fused and pooled over blocks. Each element of the generating vector is processed in parallel through this multiblock structure (aligning with its corresponding element in the opposite generating vector). These vectors are augmented by adding residuals constructed by passing them through a 1d convolution and GeLU activation, and the maximum over elements gives a feature vector for each generating vector as input for a prediction layer, answering the question whether two generating vectors lie in the same topological equivalence class. 
    }
    \label{fig:twin}
\end{figure}
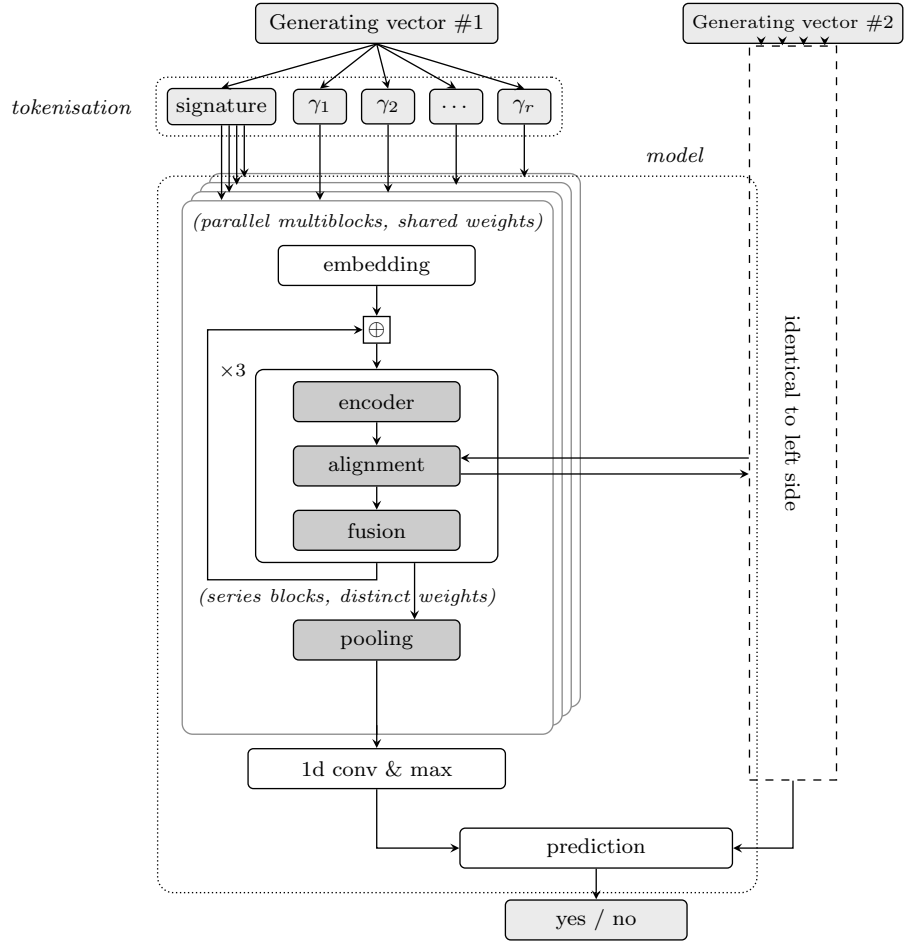

We approach the problem as one of supervised binary classification: given a pair of generating vectors sharing a signature, predict whether they lie in
the same topological equivalence class. This approach, in contrast to one which assigns a generating vector to an entry in a list of known equivalence classes, was used because a fixed-class classifier cannot by construction generalise to equivalence classes absent from its training data.
\begin{remark}
    An early proof-of-concept design as described above, one which assigned a generating vector to a fixed list of equivalence classes, performed well within its class list but was still ultimately abandoned for generalisability reasons.
\end{remark}
All models were implemented in PyTorch \cite{Paszke2019}, adapting the existing open-source implementation of the RE2 text-matching architecture \cite{Yang2019}; our final model architecture is visualised in Figure \ref{fig:twin}. Our modifications are
confined to the input side of the network (tokenisation and embedding), so
nothing hand-engineered about the group theory enters the model except
through the tokens themselves; this affords us the benefit of being able to compare the model accuracy on our task of learning equivalence classes against its accuracy on the task of sentiment analysis for which it was designed. While we shall not repeat the description of the architecture provided in \cite{Yang2019}, we will provide comments where appropriate. 

Training was undertaken on the Aire HPC system at the University of Leeds, each run occupying a single NVIDIA L40S 48GB (PCIe) GPU. 
Models were trained with the Adam optimiser under a cross-entropy loss, with learning rate $2 \times 10^{-3}$ decaying exponentially ($\times 0.94$ per $25{,}600$ training samples, floored at $10^{-4}$) after a $20{,}480$-sample warm-up, batch size 8 with 4-step gradient accumulation (an effective batch of 32), and gradient clipping at 3. Model performance was assessed using a standard train-validation-test split; throughout training, the model was evaluated periodically on the validation data and the highest-scoring checkpoint was retained, with the subsequent accuracy on the test data reported. A run terminated when its evaluation performance failed to improve for $256{,}000$ consecutive samples, or at a 24-hour wall limit.

As training data we computed all generating vectors for all group actions on Riemann surfaces of genus $2 \leq g \leq 48$ with quotient genus $g_0=0$, and found their topological equivalence classes. Data of all possible pairs of groups and signatures which act on curves of genus $2 \leq g \leq 48$ were provided by Jen Paulhus \cite{Paulhus2015}.\footnote{The link to the data in the cited paper is now outdated: the correct link at time of writing is \url{www.jenpaulhus.com/research/monodromy.html}} From this data, code adapted from \cite{Behn2023} was run in SageMath \cite{sagemath} to find all generating vectors and partition them into their equivalence classes. This resulted in 2{,}996 distinct groups, 27{,}122 signatures for these groups, 176{,}728 topological equivalence classes, and 960{,}201{,}074  generating vectors.
From this we generated non-overlapping train, validation and test datasets of pairs of generating vectors with the same signature via random sampling (uniformly over signatures, with probability $0.5$ that a pair were in the the same equivalence class), the train dataset consisting of $N=300,000$ pairs, the validation dataset consisting of $10,000$ pairs. The test dataset is a class-uniform census: for each equivalence class we sample four same-class pairs and up to four different-class pairs sharing its signature, giving $86,760$ pairs. In these datasets signatures of length $r \leq 98$ were observed (the maximum corresponding to $C_2$ acting as the hyperelliptic involution on a curve of genus $g=48$), and we also performed an evaluation limiting to signatures of length $r \leq 4$ and $r \leq 16$ to see how performance varied. Results are given in Table \ref{tab:results-embedding}.

\begin{remark}
In optimising the code for computing all the generating vectors and their topological equivalence classes, the main improvement came from precomputing many properties for each finite group $G$. In particular, we generated an ordered list of elements, and thus turned group multiplication, inverse, the order of a group element, and the action of generators of $\mathrm{Aut}(G)$ into integer arrays. This means that a single call of the computer algebra system Sage (and implicitly GAP) was required for each group, and then computing every generating vector and its orbit for each signature associate to that group could be done entirely using Python arrays. 
\end{remark}


Each element of the generating vector is then tokenised following the approach of \cite{Petschack2025}: an element of the symmetric group $S_n$ is written as a product of at most $n-1$ transpositions by decomposing $l$-cycles as 
\[
(a_1, \dots, a_l) = (a_1, a_2)(a_2, a_3) \dots (a_{l-1}, a_l),
\]
transpositions in $S_n$ are given by integers via the map 
\[
s_{ij} := (i, j) \mapsto  (i-1) + n(j-1) \in \lbrace 0, \dots, n^2-1 \rbrace, 
\]
and so the tokenisation of the element $s_{i_1 j_1} \cdots s_{i_l j_l}\in S_n$ is a sequence of tokens $t_1, \dots, t_l$, with each $t_k \in \lbrace{0, \dots, n^2-1}\rbrace$. In generating our dataset every group was given by a faithful permutation representation of degree at most $n=500$, so tokenising via the symmetric group is well defined. In practice we also add a token for the order of the element, and then append pad tokens to fix the sequence length for each group element to be the same. 

\begin{remark}
    Abstract groups were given a faithful permutation representation via the GAP commands \texttt{IsomorphismPermGroup} and \texttt{SmallerDegreePermutation\allowbreak 
    Representation}. The upper bound of $n=500$ was not strict, a lower value could be found and so one could use fewer distinct tokens to tokenise the element of generating vectors.
\end{remark}

\begin{remark}
    Transpositions were chosen in \cite{Petschack2025} because they are known to generate the symmetric group, and the maximum word length for an element of the symmetric group as written in the alphabet of transpositions is known in advance. Not all the possible tokens in $\lbrace 0, \dots, n^2-1\rbrace$ will appear, for example not $i-1$ for $i=1,\dots, n$ which corresponds to $s_{i1}$.

    One can in principle tokenise elements of the symmetric group in an analogous way provided any generating set for which one knows the maximum possible word length (and can solve the word problem efficiently). As a step in experimentation we attempted using a single transposition and 3-cycles as the generating set, but this led to too many tokens and so was abandoned as a possible tokenising step. One may ask the question as to whether different generating sets can be chosen leading to a smaller number of tokens. 
\end{remark}

The tokenisation process produces sequences of varying token length for each element of a generating vector. Given generating vectors $(\gamma_1^{(j)}, \dots, \gamma_r^{(j)})$ ($j=1, 2$) being compared the model generates feature vectors for $\gamma_i^{(1)}$ and $\gamma_i^{(2)}$ using a copy of the embedding-encoder-fusion-pooling  process implemented in the text-matching architecture of \cite{Yang2019}. 
\begin{remark}
    The model in \cite{Yang2019} consists of a series of repeated blocks which process the input. As in the original paper we experimented with varying the number of blocks to be 1, 2, 3 and 4, seeing only marginal impact upon performance. As such, we kept the default value of 3 blocks.
\end{remark}
From these feature vectors $v_i$ for each element $\gamma_i$, a feature vector for the whole generating vector is found as follows: $(v_1, \dots, v_r)$ is passed through a 1d convolution layer followed by GeLU activation to give a residual $c_i$, and then an elementwise max is applied to $v_i+c_i$ across each element of the generating vector. The prediction layer from \cite{Yang2019} in symmetric configuration is used to get an answer from the two feature vectors corresponding to the input pair. 

\begin{remark}
    The use of the 1d convolution before applying the maximum over element in a generating vector ensures that the model in not invariant under permuting the elements of the generating vector, as in this latter scenario model errors are guaranteed from configurations such as in Example \ref{ex: bad gv}. In development we experimented with positional embedding features but found that this performed worse.

    Conversely, using the symmetric prediction layer means that the model prediction for whether two generating vectors are in the same equivalence class does not depend on the order they are provided to the model as one would desire. 
\end{remark}

We compared two methods of padding generating vectors so that a variable length $r$ was accepted,
\begin{enumerate}
    \item padding the generating vector with copies of the identity element subsequently tokenised, effectively fixing the generating vector length to be some maximum $r$, and 
    \item padding the generating vector with PAD tokens subsequently ignored using an attention mask in the softmax and pooling layer.
\end{enumerate}
The results of this padding comparison are given in Table \ref{tab:results-embedding}.

Tokens enter the network through a trainable embedding table, initialised in two different ways,
\begin{enumerate}
    \item analytically, from the $24 \times 24$ matrix representation induced from the permutation representation of $S_{24} \leq S_{500}$: the token of a transposition is initialised to its flattened permutation matrix, giving embedding vectors of dimension $576$, and  
\item randomly initialised with dimension $576$, chosen to allow for comparison with the matrix representation, where any difference in performance between the two initialisations may be attributed structure rather than size.
\end{enumerate}
The results of this comparison are given in Table~\ref{tab:results-embedding}.

Fixing the use of the padding mask as the method for handling variable $r$ and initialising with the random embedding, we wished to explore how the model would perform generalising to generating vectors considerably different to the training data. To do this, we trained the model only on groups in the training data of order $\abs{G} \leq 64$, validating and testing on groups of order $\abs{G} > 64$. The value of 64 was chosen to strike a balance between seeing sufficiently many smaller groups while leaving a big enough variety of larger groups to test on; the number of distinct finite groups of order $n$ peaks considerably at $n=64$, and approximately 80\% of signatures in our dataset occur for $\abs{G} \leq 64$, thus informing this choice. With this configuration, again the validation and test datasets contained $10,000$ and $86,760$ pairs respectively, but we varied the size of the training dataset $N$ to investigate how this limited performance. The results of this are given in Table \ref{tab:results-scaling}.

To better understand any errors in classification that occurred in the latter test stage (including only pairs generating groups of order $\abs{G}>64$) we labelled the pairs with data on the groups they generate from the LMFDB \cite{lmfdb}, and then compared model performance (when trained on $N=10^6$ pairs generating groups of order $\abs{G} \leq 64$) on subclasses of pairs by plotting Receiver Operating Characteristic (ROC) curves (Figures \ref{fig:ROC1} and \ref{fig:ROC2}) with the corresponding Area-Under-Curve (AUC) values. In particular we focused on the following subclasses of generating vector pairs:
\begin{itemize}
    \item cross-group, where the generating vectors generate different groups, 
    \item same-group non-abelian, where the generating vectors generate the same non-abelian group, 
    \item same-group cyclic, where the generating vectors generate the same cyclic group, 
    \item same-group other abelian, where the generating vectors generate the same non-cyclic abelian group, 
    \item same-group $p$-group, where the generating vectors generate the same $p$-group.
\end{itemize}
Moreover, to investigate patterns in the feature vectors found by the model we performed dimensional reduction using t-distributed Stochastic Neighbour Embedding (t-SNE) plots of each pair's feature vectors, in particular using the absolute elementwise difference and product  (Figures \ref{fig:tSNE_diff} and \ref{fig:tSNE_prod}), as these are the symmetric combinations of the features for each generating vector which enter the prediction layer.

\section{Results}

\subsection{Classification accuracy}

Table ~\ref{tab:results-embedding} reports the test accuracy of each training configuration when not limiting $\abs{G}$ under the aforementioned train-validation-test split, organised by padding method, maximum $r$ and embedding choice. The best performance reported was an accuracy of $85.8\%$.
We see that the choice of embedding initialisation only made a meaningful difference ({$\sim 2\%$}) when using identity padding or for larger $r$ where the padding mask was used, but that uniformly the random embedding gave better performance.
We were unable to test the accuracy of the model with identity padding for larger maximum values of $r$ due to memory constraints.
Comparing $r \leq 16 $ and $r \leq 98 $ with the padding mask we saw that there was no significant difference in accuracy.

\begin{table}
\centering
\begin{tabular}{lcc}
\hline
Signature Length & \multicolumn{2}{c}{Accuracy} \\
\hline & \underline{Random Embedding} & \underline{Matrix Embedding} \\

$r \leq 4^\ast$
& $0.702$  
& $0.681$ \\
$r \leq 4$
& $0.698$  
& $0.693$ \\
$r \leq 16$
& $0.850$ 
& $0.828$\\ 
$r \leq 98$
& $\bm{0.858}$ 
& $0.825$ \\ 
\hline
\end{tabular}
\caption{Model accuracy on $86,760$ test pairs when trained on $300,000$ pairs with $\abs{G}$ unrestricted and varying range for $r$. ${}^\ast$ indicates that identity padding was used. }
\label{tab:results-embedding}
\end{table}

Table \ref{tab:results-scaling} shows the performance of the model restricted to training on a varying number of generating vectors from groups of order $\abs{G} \leq 64$ when evaluated on a test dataset consisting of pairs of generating vectors from groups with $\abs{G} > 64$. For these experiments the embedding initialisation was fixed to be random, and the padding method used an attention mask. Model performance is seen to plateau after reaching approximately $0.795$.

\begin{table}
\centering
\begin{tabular}{lccc}
\hline
Training Pairs & Accuracy\\
\hline
$1\times 10^5$    & $0.743$ \\ 
$2\times 10^5$    & $0.767$ \\ 
$3\times 10^5$    & $0.789$ \\ 
$5\times 10^5$    & $0.781$ \\ 
$8\times 10^5$    & $0.792$ \\ 
$1\times 10^6$    & $\bm{0.795}$ \\ 
$1.3\times 10^6$  & $0.795$ \\ 
$2\times 10^6$    & $0.792$ \\ 
\hline
\end{tabular}
\caption{Accuracy against training set size
(train $|G| \leq 64$, test $|G| > 64$)}
\label{tab:results-scaling}
\end{table}

\subsection{Error analysis}

We analyse the errors of the random-embedding model on the test set. Of the $86,760$ pairs, the model misclassified $17,786$. The error rates on specific subclasses are given in Table \ref{tab:error by class}. Of the subclasses considered the model performs best on the same-group $p$-group subclass, closely followed by the cross-group subclass. The model performs similarly well on the same-group non-abelian class, and conversely performs poorly for all same-group abelian pairs both cyclic and otherwise.

\begin{table}
\centering
\begin{tabular}{ccc}
\hline
Subclass & Accuracy & Number of Pairs \\
\hline
all pairs               & $0.795$      & $86,760$ \\ 
cross-group             & $0.823$      & $33,654$ \\ 
same-group non-abelian  & $0.813$      & $44,861$ \\ 
same-group cyclic       & $0.594$      & $7,184$  \\ 
same-group other abelian& $0.557$      & $1,061$  \\ 
same-group $p$-group    & $\bm{0.829}$ & $5,146$  \\ 
\hline
\end{tabular}
\caption{Model accuracy on subclasses of generating vector pairs}
\label{tab:error by class}
\end{table}

\begin{figure}
    \centering  
    \includegraphics[width=0.95\textwidth]{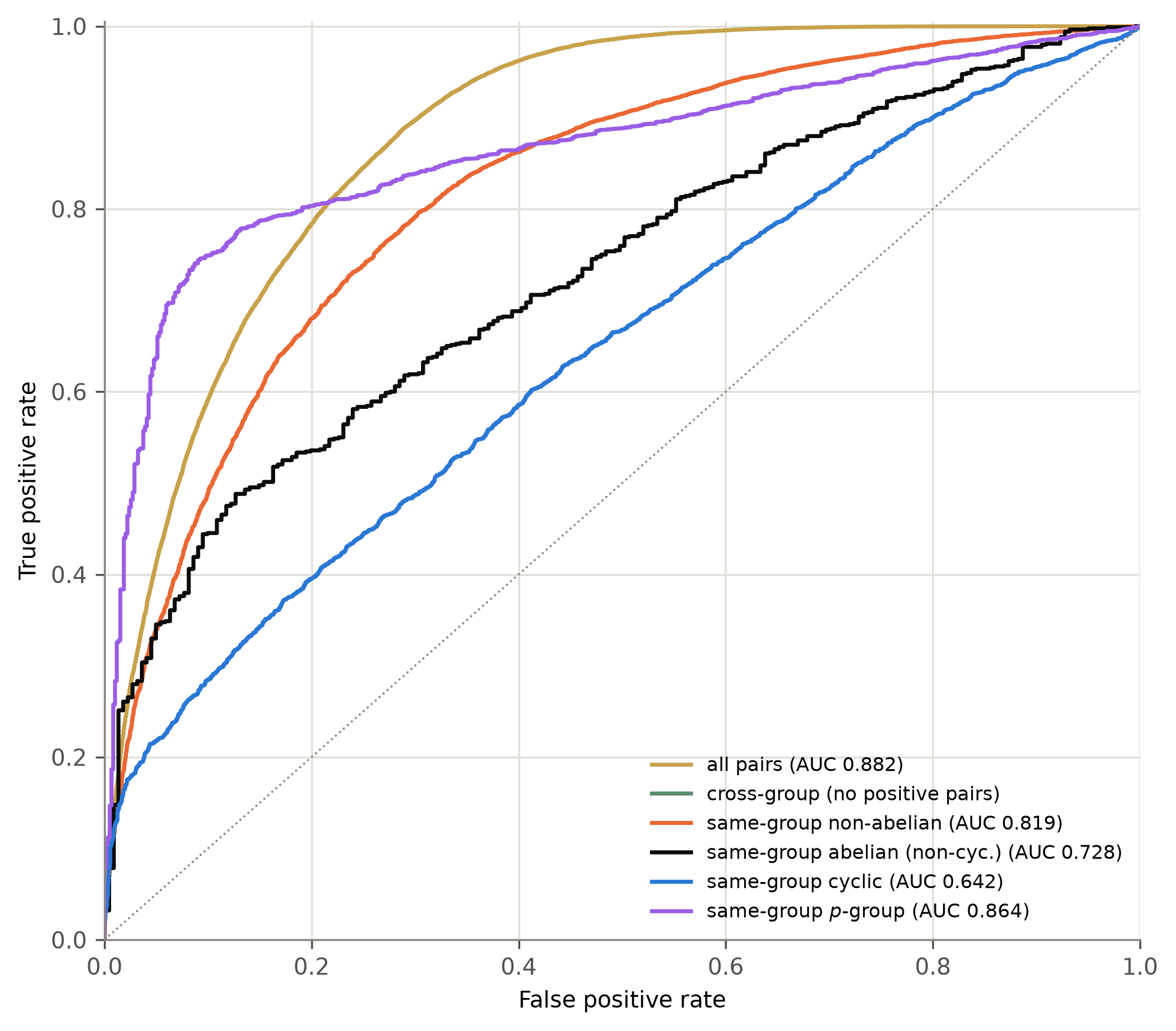}
    \caption{ROC curves for the model classifying pairs as in the same equivalence class. Subclasses are labelled by the AUC.  
    }
    \label{fig:ROC1}
\end{figure}

Figures \ref{fig:ROC1} varies the confidence threshold in the output of the model for classification of pairs as in the same equivalence class from $1$ to $0$, and compares the true positive rate against the false negative rate. Overall the model performs reasonably well as a classifier on all pairs with an AUC of $0.882$; 
this is noticeably higher than the performance on most same-group subclasses considered, indicating that the false-positive rate is largely suppressed by cross-group pairs which are accurately identified as in different equivalence classes. The exception is the same-group $p$-group subclass, where the high-confidence positive likelihood ratio exceeds that on all pairs.

\begin{figure}
    \centering  
    \includegraphics[width=0.95\textwidth]{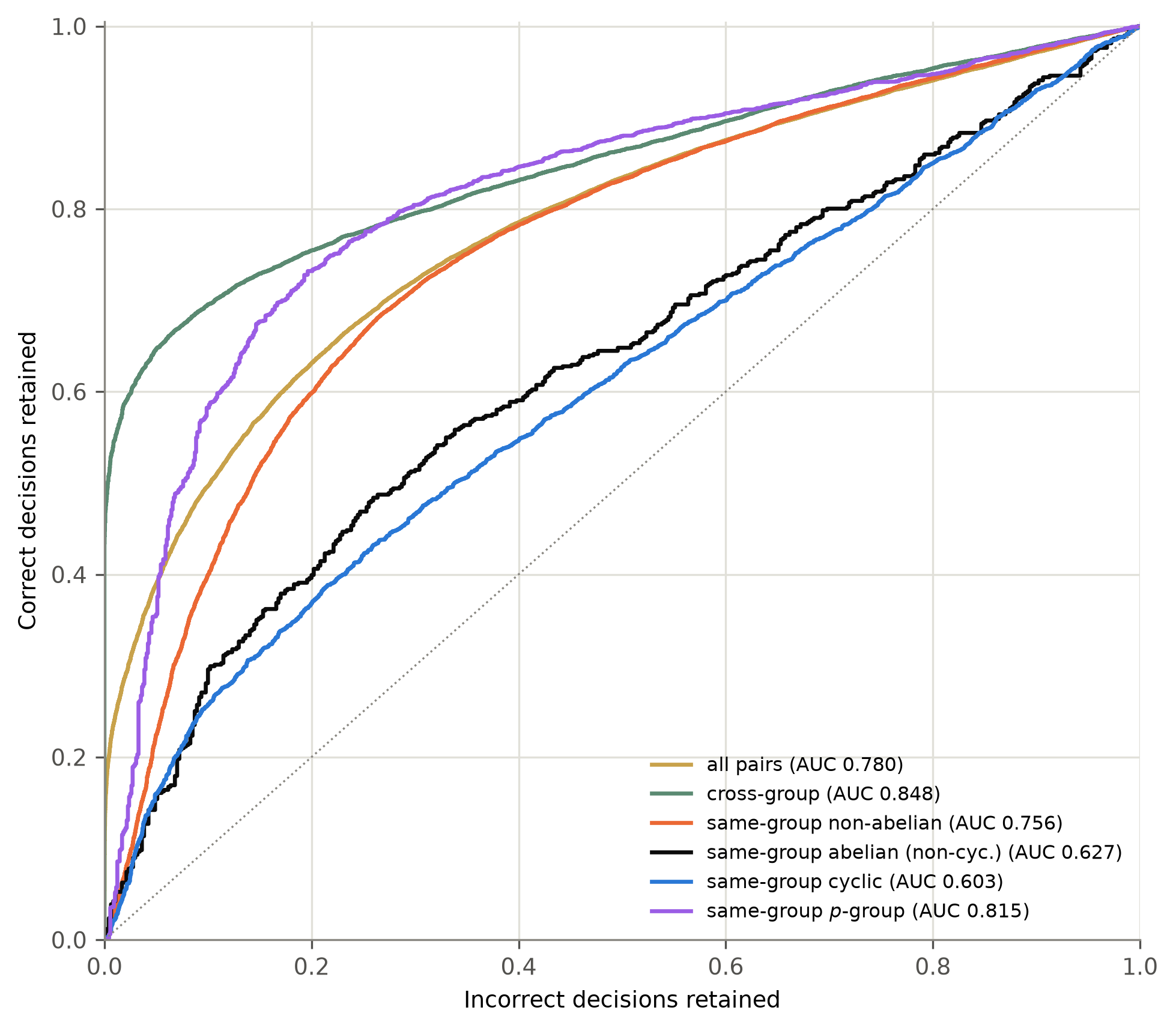}
    \caption{ROC curves for the model classifying pairs as the confidence threshold is varied. Subclasses are labelled by the AUC.
    }
    \label{fig:ROC2}
\end{figure}

Because there is no `true positive' for pairs in the cross-group subclass when classifying whether pairs are in the same equivalence class, Figure \ref{fig:ROC2} is a modification of the ROC curve varying the confidence threshold for classification of pairs from $1$ to $0.5$ (the default value) and comparing how many correct and incorrect decisions are retained; for this curve we may consider the cross-group subclass. One must be careful interpreting the AUC for these curves: even a value near 1 could correspond to a subclass of pairs on which the model performed very poorly if there were many incorrect decisions with confidence near $0.5$. We nevertheless see that the model confidently distinguishes a large proportion of cross-group and same-group $p$-group pairs. 

\begin{figure}
    \centering  
    \includegraphics[width=0.95\textwidth]{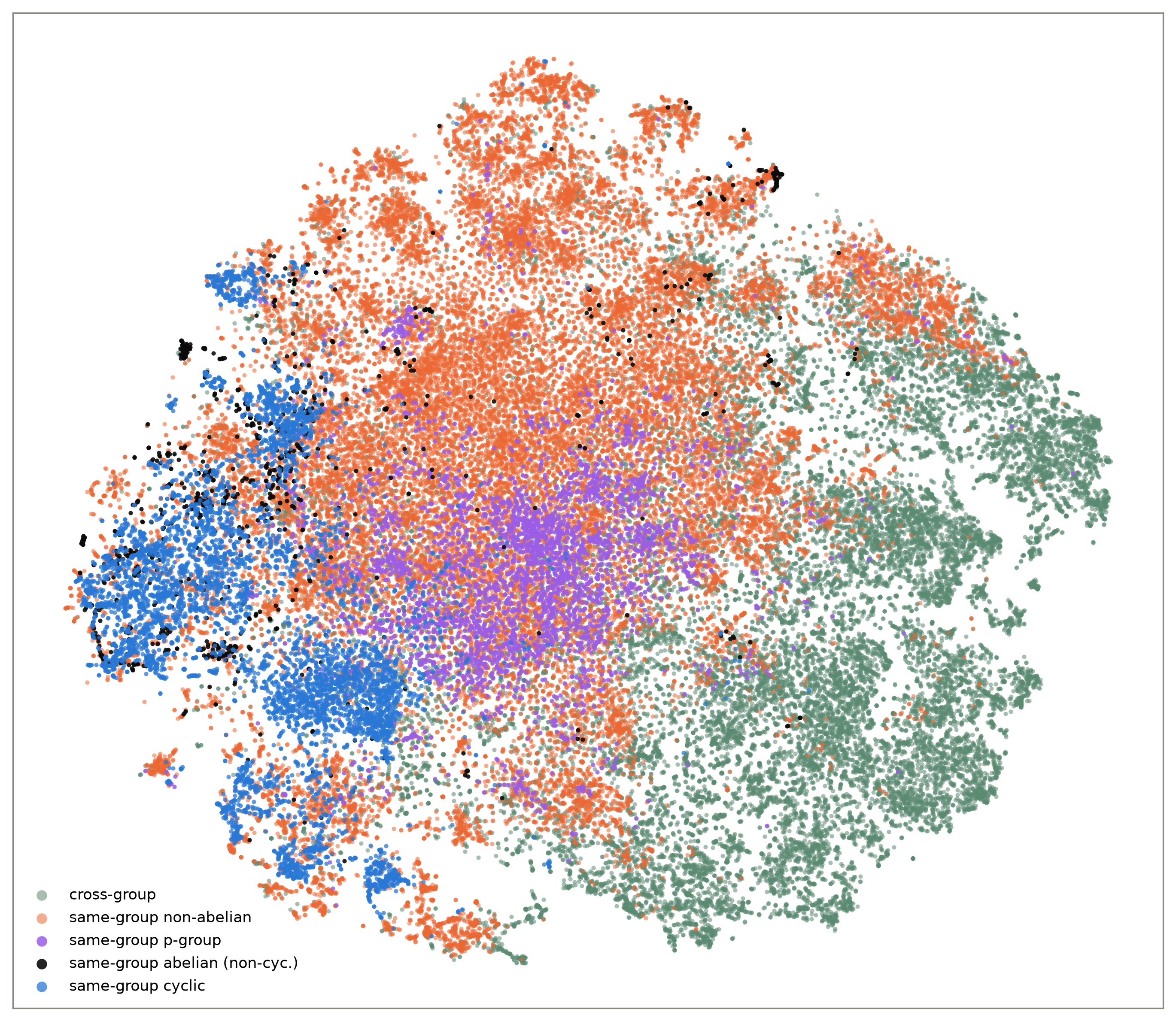}
    \caption{Dimensional reduction of the absolute elementwise difference in the generating vector features of a pair, stratified by subclass
    }
    \label{fig:tSNE_diff}
\end{figure}

\begin{figure}
    \centering  
    \includegraphics[width=0.95\textwidth]{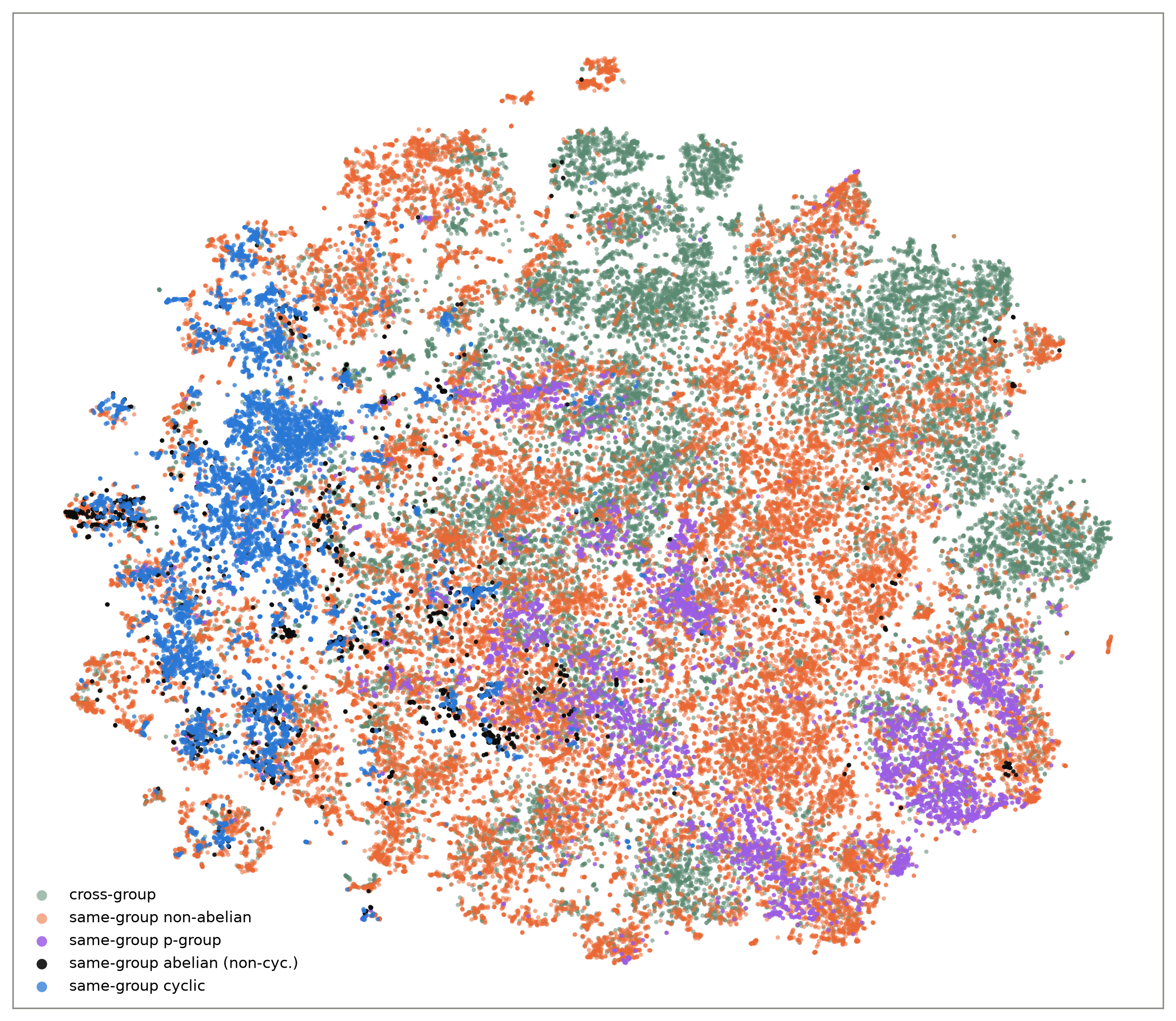}
    \caption{Dimensional reduction of the elementwise product of the generating vector features of a pair, stratified by subclass
    }
    \label{fig:tSNE_prod}
\end{figure}

Figures \ref{fig:tSNE_diff} and \ref{fig:tSNE_prod} partially explain the model's respective high and low performance on the cross-group and same-group other abelian subclasses: the former are well separated by the model in feature space, whereas the latter are not tightly clustered. This perspective further provides a plausible explanation for the marginally higher performance on the same-group $p$-group subclass than on the same-group non-abelian subclass which contains it, but is unable to explain why the performance on same-group cyclic pairs was poor, suggesting that the prediction layer is indeed identifying local relational properties between the two features. 

\section{Discussion}

This paper illustrates how simple deep learning approaches developed in separate contexts may be applied to achieve high accuracy predictions of sophisticated group structure. In particular, a deep learning model developed for sentiment analysis of text has been adapted to learn topological equivalence of generating vectors with an accuracy of $85.8\%$. 
Notably this accuracy is comparable with that of the model when evaluated on a text-matching task given similarly size train data. The original RE2 architecture has accuracy $89.2 \pm 0.2$ on the task of identifying whether, given $400,000$ training pairs of questions from Quora, one is a paraphrase of the other \cite{Yang2019}. This is the closest analogue to this paper's task in the various tests of \cite{Yang2019}: a binary decision on whether two token sequences are equivalent. Our accuracies when training on $300,000$ pairs with $r \leq 98$ and $|G|$ unrestricted (see Table~\ref{tab:results-embedding}) are directly comparable to the original architecture's performance. The comparison is not exact: there are differences in class balance, embedding dimension, vocabulary size and sequence structure. Nonetheless, there is a striking observation that the architecture learns topological equivalence of generating vectors comparably well to how it learns the semantic equivalence of English sentences, the task it was designed for. This suggests that the model has been able to identify some level of semantic structure in topological equivalence classes.

Unlike the application in \cite{Yang2019} which used the 300-dimensional GloVe word vectors, we predominantly did not use a tailored pre-trained embedding as this did not improve model performance. We hypothesise that this is because the matrix embedding we experimented with, implicitly providing the data of a 24-dimensional representation of $S_{500}$, is insufficiently helpful for the task of identifying topological equivalence classes. Indeed such a representation cannot be faithful as the minimal dimension for a faithful representation of $S_n$ is $n-1$: this would correspond in our case to embedding vectors of dimension $499^2=249001$ which would have been impracticably large. Our results do not rule out the possibility that another tailored embedding would provide improvement, and future work should explore which tokenisations and embeddings of group elements are best tailored to learning finite group theory.

We saw that, though the model performance did drop when trained only on groups of order $\abs{G} < 64$, it was still able to generalise reasonably well. This bolsters the conclusion that some underlying structure of topological equivalence was learnt, though one must caveat this with the fact that, though a given group is not present in both the train and test data, an element of that group may be present in a generating vector in both. Though accuracy improved with the size of training data available as is natural it eventually plateaued, indicating that there are still aspects of topological equivalence which were missed. 
It is possible, if apparently unlikely, that with a much larger training dataset performance would begin to increase again and cross a $99\%$ accuracy threshold, so-called `grokking' \cite{Power2022}. As observed in \cite{Power2022}, this behaviour requires the amount of training data to be approximately at least a third of all possible data: a rough estimate suggests that from our dataset one could generate $O(10^{13})$ pairs of generating vectors with the same signature, many orders of magnitude more than were in the largest train dataset.

Data mining tools have demonstrated that a key factor underlying the models strength even when trained only on pairs generating groups of order $\abs{G} \leq 64$ is that the model is confidently able to distinguish those pairs which generate different groups, even though this information was not a priori provided to the model. Dimensional reduction indicates that the feature vectors of such pairs are well separated from other same-group pairs. What error analysis has highlighted is that the model performs well on the same-group $p$-group subclass, and that this is not immediately attributable to a separation of feature vectors. We still lack an explanation for why this behaviour occurs. The model's relative underperformance on abelian groups (both cyclic and non-cyclic) would suggest that the model has not developed the ability to identify when the ordering of elements becomes immaterial, even though the feature vectors for the same-group cyclic subclass are well separated.

It is worthwhile commenting further on the decision to use a deep learning approach in contrast to a machine learning approach, here the distinction being that no feature engineering was undertaken to extract significant data from the input generating vectors beforehand on which to train, instead resorting to tokenisation. In this aspect, our paper is more closely aligned with the approach of \cite{Petschack2025} than with that of \cite{Braden2025}. In doing so we increase the flexibility of the model at the cost of losing interpretability, and we have been unable to extract any new theorems from this study as of yet. That being said, our error analysis indicates that the model best learnt the topological equivalence classes of $p$-groups in the dataset, which may indicate that this class of groups is a natural candidate for future mathematical study.
Moreover, deep learning is typically more computationally intensive than machine learning, and runtime has prevented our ability to evaluate the performance of the model when trained on more than a tiny fraction of the total possible data.

By design we did not attempt to force invariance of the model under the action of $\Aut(G) \times \mathbb{B}_r$ on either generating vector input (appropriately understood), hoping that our model would learn this structure. Employing this, future work could investigate the application of a multiple-testing framework to determine if two generating vectors are in the same equivalence class, comparing multiple pairs of generating vectors from either class at a time.


\bibliography{refs}
\addcontentsline{toc}{chapter}{Bibliography}
\bibliographystyle{amsalpha}

\end{document}

%% file: twin_network_tikz_paper_v2.tex
\begingroup
\begin{tikzpicture}[
    font=\footnotesize,
    line width=0.5pt,
    >=stealth,
    data/.style={draw, fill=black!8, rounded corners=2pt,
                 minimum height=15pt, minimum width=3.2cm},
    tok/.style={draw, fill=black!8, rounded corners=2pt,
                minimum height=13pt, minimum width=7mm},
    layer/.style={draw, fill=white, rounded corners=2pt,
                  minimum height=15pt, minimum width=2.6cm},
    stage/.style={draw, fill=black!20, rounded corners=2pt,
                  minimum height=15pt, minimum width=2.2cm},
    opbox/.style={draw, fill=white, minimum size=10pt, inner sep=1pt},
    lab/.style={font=\footnotesize\itshape},
    shadowbox/.style={draw=black!45, fill=white, rounded corners=4pt},
  ]

  \node[data] (gv) at (0, 0) {Generating vector \#1};
  \node[tok, minimum width=1.4cm] (t0) at (-2.05, -1.1) {signature};
  \node[tok] (t1) at (-0.75, -1.1) {$\gamma_1$};
  \node[tok] (t2) at ( 0.15, -1.1) {$\gamma_2$};
  \node[tok] (t3) at ( 1.05, -1.1) {$\dots$};
  \node[tok] (t4) at ( 1.95, -1.1) {$\gamma_r$};
  \node[draw, densely dotted, rounded corners=4pt,
        fit=(t0)(t4), inner sep=4pt] (tokbox) {};
  \node[lab, anchor=east] at ([xshift=-2mm]tokbox.west) {tokenisation};
  \foreach \t in {t0,t1,t2,t3,t4}
    \draw[->] (gv.south) -- (\t.north);

  \node[shadowbox, minimum width=4.9cm, minimum height=7.05cm]
    at (0.24, -5.51) {};
  \node[shadowbox, minimum width=4.9cm, minimum height=7.05cm]
    at (0.12, -5.63) {};
  \node[shadowbox, minimum width=4.9cm, minimum height=7.05cm]
    (pergen) at (0, -5.75) {};
    \node[shadowbox, minimum width=4.9cm, minimum height=7.05cm]
    (pergen) at (-0.12, -5.87) {};
  \node[lab, font=\scriptsize\itshape, anchor=north west]
    at ([xshift=0mm, yshift=-0.5mm]pergen.north west)
    {(parallel multiblocks, shared weights)};

  \draw[->] (t0.south) -- ([yshift=0.0cm]t0.south |- pergen.north);
  \draw[->] ([xshift=0.1cm]t0.south) -- ([yshift=0.1cm, xshift=0.1cm]t0.south |- pergen.north);
  \draw[->] ([xshift=0.2cm]t0.south) -- ([yshift=0.2cm, xshift=0.2cm]t0.south |- pergen.north);
  \draw[->] ([xshift=0.3cm]t0.south) -- ([yshift=0.3cm, xshift=0.3cm]t0.south |- pergen.north);
  \draw[->] (t1.south) -- ([yshift=0.0cm]t1.south |- pergen.north);
  \draw[->] (t2.south) -- ([yshift=0.1cm]t2.south |- pergen.north);
  \draw[->] (t3.south) -- ([yshift=0.2cm]t3.south |- pergen.north);
  \draw[->] (t4.south) -- ([yshift=0.3cm]t4.south |- pergen.north);

  \node[layer] (emb) at (0, -3.2) {embedding};
  \node[opbox] (oplus) at (0, -4.05) {$\oplus$};
  \draw[->] (emb.south) -- (oplus.north);

  \node[draw, rounded corners=3pt, minimum width=3.2cm, minimum height=2.55cm]
    (blk) at (0, -5.85) {};
  \node[lab, anchor=north east, font=\scriptsize\itshape]
    at ([yshift=-2mm, xshift=1mm]blk.south east) {(series blocks, distinct weights)};
  \node[stage] (enc) at (0, -5.0) {encoder};
  \node[stage] (ali) at (0, -5.85) {alignment};
  \node[stage] (fus) at (0, -6.7) {fusion};
  \draw[->] (oplus.south) -- (blk.north);
  \draw[->] (enc.south) -- (ali.north);
  \draw[->] (ali.south) -- (fus.north);
  \draw[->] (blk.south) -- ++(0, -0.22) -| ([xshift=-2.05cm]oplus.west)
    |- (oplus.west);
  \node[lab, font=\scriptsize\itshape, anchor=east] at (-1.6, -4.6) {$\times 3$};

  \node[stage] (pool) at (0, -8.15) {pooling};
  \draw[->] ([xshift=0.5cm]blk.south) -- ([xshift=0.5cm]blk.south |- pool.north);

  \node[layer, minimum width=3.4cm] (mx) at (0, -9.85) {1d conv \& max};
  \draw (pool.south) -- (pool.south |- pergen.south);
  \draw[->] (0, -9.275) -- (mx.north);
  \node[layer, minimum width=3.6cm] (pred) at (2.9, -10.9) {prediction};
  \draw[->] (mx.south) |- (pred.west);

  \node[draw, dashed, minimum width=1.15cm, minimum height=9.7cm]
    (twin) at (5.5, -5.15) {};
  \node[rotate=-90, font=\footnotesize] at (twin.center)
    {identical to left side};
  \node[data, minimum width=2.9cm, font=\scriptsize] (gv2) at (5.5, 0)
    {Generating vector \#2};
  \foreach \dx/\st in {-0.42/solid, -0.14/solid, 0.14/densely dotted, 0.42/solid}
    \draw[->, \st] ([xshift=\dx cm]gv2.south)
      -- ([xshift=\dx cm]gv2.south |- twin.north);
  \draw[<-] ([yshift=3pt]ali.east) -- ([yshift=3pt]twin.west |- ali.east);
  \draw[->] ([yshift=-3pt]ali.east) -- ([yshift=-3pt]twin.west |- ali.east);
  \draw[->] (twin.south) |- (pred.east);

  \node[data, minimum width=2.4cm] (yn) at (2.9, -11.85) {yes / no};
  \draw[->] (pred.south) -- (yn.north);

  \node[draw, densely dotted, rounded corners=6pt, inner sep=9pt,
        fit=(pergen)(mx)(pred)] (model) {};
  \node[lab, anchor=south east]
    at ([yshift=0.5mm, xshift=-6mm]model.north east) {model};
\end{tikzpicture}
\endgroup